\documentclass{elsart}
\usepackage{graphicx}
\usepackage{epsfig}
\usepackage{amssymb}

\begin{document}

\begin{frontmatter}
\thanks[label1]{Thanks to our visit to Osaka City University in summer
2004 sponsored by COE program of Prof Akio Kawauchi where we learnt
about the quasi toric braids and the Manturov's theorem.}

\title{Polynomial Representation for Long Knots}

\author{Rama Mishra$^{\dag}$ and M. Prabhakar$^{\ast}$}

\address{
$\dag$ Indian Institute of Science Education and Research, Pune,
India.} \ead{r.mishra@iiserpune.ac.in}
\address{$\ast$ Department of Mathematics,
Indian Institute of Technology Guwahati, Guwahati - 781 039, India}
\ead{prabhakar@iitg.ernet.in}

\begin{abstract}
We discuss the polynomial representation for long knots and
elaborate on how to obtain them with a bound on degrees of the
defining polynomials, for any knot-type.
\end{abstract}


\end{frontmatter}

\section{Introduction}
This paper is aimed to be a survey article on the polynomial
representation of knots and includes all the results proved by the
authors. A. Durfee and D. Oshea \cite{ad} wrote a similar paper in
2006. Our paper provides a more constructive approach and uses
recent knot theoretic results. Polynomial Representation for long
knots were first shown by Shastri \cite{ars}. Shastri proved that
for every knot-type $K \ (\mathbb{R} \hookrightarrow \mathbb{R}^3)$
there exist real polynomials $f(t), g(t)$ and $h(t)$ such that the
map $t \mapsto (f(t), g(t), h(t))$ from $\mathbb{R}$ to
$\mathbb{R}^3$ represents $K$ and in fact the above map defines an
embedding of $\mathbb{C}$ in $\mathbb{C}^3.$ Shastri's motivation
for proving this theorem was perhaps to find a non-rectifiable
polynomial embedding of complex affine line in complex affine space
which could prove one of the famous conjecture of Abhayankar
\cite{ss1}. However, this conjecture is still open. A knot
represented by a polynomial embedding is refered as a polynomial
knot. Similar to the case of Harmonic knots \cite{tr} and Holonomic
knots \cite{ek} it has been proved that the space of all knots
\cite{am}, up to ambient isotopy, can be replaced by the space of
all polynomial knots, up to polynomial isotopy \cite{rs}. However
polynomial parametrization of knots not only represents the knot
type, it also can represent a specific knot diagram. For instance,
it has been shown that we can capture the symmetric behavior of
knots such as strong invertibility and strong negative
amphicheirality in their polynomial representation \cite{rm3}.
Polynomial knots are also important from the point of view of
computability as  we can easily obtain three dimensional graphs of
polynomial knots using {\it Mathematica} or {\it Maple} as
illustrated in section 6.

 In Shastri's theorem the existence of a polynomial representation
for a given knot type $K$ is shown by using Weierstrass'
approximation. Thus, estimating the degrees of the polynomials is
not clear. Suppose $t \mapsto (f(t), g(t), h(t))$ is a polynomial
representation of a knot type $K$ and $deg(f(t)) = l, \ deg(g(t)) =
m$ and $deg(h(t)) = n$, then we say that $(l, m, n)$ is a {\it
degree sequence} of $K$. We define $(l, m, n)$ to be the {\it
minimal degree sequence} of $K$ if $(l, m, n)$ is minimal amongst
all the degree sequences of $K$ with respect to the usual
lexicographic ordering of $\mathbb{N}^3$. Note that a degree
sequence of $K$ need not be unique. Also, with a given triple
$(l,m,n)$ of positive integers, the space of polynomial knots with
$(l,m,n)$ as a degree sequence is a {\it semi algebraic set}
\cite{rb}, and thus only finitely many knot-types can be realized
with a given degree sequence. Hence, for a given knot-type
estimating a degree sequence and eventually the minimal degree
sequence is an important aspect for polynomial knots. In our earlier
papers (\cite{am},\cite{rm1},\cite{rm2}, \cite{p1},\cite{p2},
\cite{p3}) we have estimated a degree sequence and the minimal
degree sequence for torus knots and 2-bridge knots.

In this paper, for any given knot type we give a concrete algorithm
for writing down a polynomial representation and hence estimating a
degree sequence. For this purpose, we use a recent theory of {\it
quasi toric} braid representation of knots \cite{vom}. It is
established by Manturov that every  knot is ambient isotopic to some
knot which is obtained by making a few crossing changes in some
torus knot \cite{vom}. Thus, when we know a degree sequence for all
torus knots, we can have some estimate on a degree sequence for a
general knot also. Here we discuss this idea in detail and provide a
method  to obtain a degree sequence for a general knot type.

In section 2, we provide the basic definitions and historical
background of the subject. We  also include some known results
without the proofs. In section 3, we state and prove the main
result. We demonstrate this method and construct a polynomial
representation for the knot $8_{17}$ in section 4. Polynomial
representation of all knots up to 8 crossings and their 3d plots
taken with the help of {\it Mathematica} are included in section 5
and section 6 respectively.

\section{Some Known Results}

In this section we present the basic definitions, background and the
known results which serves as the prerequisite for the main result.
\begin{defn}
Two non-compact knots $ (\widetilde{\phi_1}: \mathbb{R}
\longrightarrow \mathbb{R}^3$) and $ (\widetilde{\phi_2}: \mathbb{R}
\longrightarrow \mathbb{R}^3)$ are said to be equivalent if their
extensions $\phi_1$ and $\phi_2$ from $S^1$ to $S^3$ are ambient
isotopic. An equivalence class of a non-compact knot is a knot-type.
\end{defn}

Once we have an embedding of $\mathbb{R}$ in $\mathbb{R}^3$ we may
wish to see if it can be extended as an embedding of $\mathbb{C}$ in
$\mathbb{C}^3$.

\begin{rem}
When we have an embedding  given by $t\mapsto (f(t),g(t),h(t))$ from
$\mathbb{R}$ to $\mathbb{R}^3,$ where $f(t),$ $g(t)$ and $h(t)$ are
three real polynomials, it defines a non-compact knot and it can
always be made into a polynomial embedding of $\mathbb{C}$ in
$\mathbb{C}^3$ by perturbing the coefficients of any one of the
polynomials.
\end{rem}

\begin{defn}
Two Polynomial embeddings $\phi_1,\:
\phi_2:\mathbb{C}\hookrightarrow \mathbb{C}^3$ are said to be {\it
algebraically equivalent} if there exists a polynomial automorphism
$F:\mathbb{C}^3\longrightarrow \mathbb{C}^3$, such that $F\circ
\phi_1=\phi_2$.
\end{defn}

\begin{defn}
A polynomial embedding $\mathbb{C}\hookrightarrow\mathbb{C}^3$ is
said to be {\it rectifiable} if it is algebraically equivalent to
the standard embedding $t\mapsto (0,0,t)$ of $\mathbb{C}$ in
$\mathbb{C}^3.$
\end{defn}

\begin{conj}(Abhayankar\cite{ss1}) There exist
{\it non-rectifiable} embeddings of $\mathbb{C}$ in $\mathbb{C}^3.$
\end{conj}

The above stated conjecture is still open.

\begin{rem}
If a polynomial embedding  $\mathbb{R}\hookrightarrow\mathbb{R}^3$
which is also an embedding of
$\mathbb{C}\hookrightarrow\mathbb{C}^3$ is rectifiable, using a
polynomial automorphism with real coefficients only, it is certainly
a trivial knot.
\end{rem}

Thus to obtain examples of non-rectifiable embeddings we must have
polynomial embeddings which represent non-trivial knots. This raises
the following question: \textit{Does every knot have a polynomial
representation?}

Fortunately we have an affirmative answer and the following results
have been proved in this line.

\begin{thm}\cite{ars}
Every knot-type ($\mathbb{R} \hookrightarrow \mathbb{R}^3$) has a
polynomial representation $t \mapsto (f(t), g(t), h(t))$ which is
also an embedding of $\mathbb{C}\hookrightarrow \mathbb{C}^3$
\end{thm}

\begin{thm}\cite{rs}
Two polynomial embeddings $\phi_0,\: \phi_1:
\mathbb{R}\hookrightarrow \mathbb{R}^3$ representing the same
knot-type are polynomially isotopic. By polynomially isotopic we
mean that there exists
$\{{P_t:\mathbb{R}\hookrightarrow\mathbb{R}^3}| \ t \in [0, 1] \},$
a one parameter family of polynomial embeddings, such that $P_0 =
\phi_0$ and $P_1 = \phi_1.$
\end{thm}

Both these theorems were proved using \textit{Weierstrass'
approximation}. Thus the nature and the degrees of the defining
polynomials cannot be estimated.

\begin{defn}
A triple $(l,m,n)\in\mathbb{N}^3$ is said to be a {\it degree
sequence} of a given knot-type $K$ if there exists $f(t),\ g(t)$ and
$h(t),$ real polynomials, of degrees $l$, $m$ and $n$ respectively,
such that the map $t\mapsto (f(t),g(t),h(t))$ is an embedding which
represents the knot-type $K$.
\end{defn}

\begin{defn}
A degree sequence $(l,m,n)\in\mathbb{N}^3$ for a given knot-type is
said to be the minimal degree sequence if it is  minimal amongst all
degree sequences for $K$ with respect to the lexicographic ordering
in $\mathbb{N}^3.$
\end{defn}

The following results have been obtained in the past.

\begin{thm}\cite{am} A torus knot of type $(2, 2n+1)$ has a degree sequence
$(3,4n,4n+1).$
\end{thm}

\begin{thm}
\cite{rm1} A torus knot of type $(p, q),  p <q$, $p>2$ has a degree
sequence $(2p-1,2q-1,2q)$.
\end{thm}

It is easy to observe that these degree sequences are not the
minimal degree sequence for torus knots. For minimal degree sequence
we have the following:

\begin{thm}
\cite{rm2} The minimal degree sequence for torus knot of type $(2,\
2n + 1)$ for $n=3m;\ 3m+1  \ \;and\; \ 3m+2$ is   $(3,\ 2n+2,\
2n+4)$; \ $(3,\ 2n+2,\ 2n+3)$ and $(3,\ 2n+3,\ 2n+4)$ respectively.
\end{thm}

\begin{thm}
\cite{p2} The minimal degree sequence for a 2-bridge knot having
minimal crossing number $N$ is given by
\begin{enumerate}
\item $(3,\ N+1,\ N+2)$ when $N \equiv 0\ (mod \ 3)$;
\item $(3,\ N+1,\ N+3)$ when $N \equiv 1\ (mod \ 3)$;
\item $(3, N+2, N+3)$ when $N \equiv 2\ (mod \ 3)$
\end{enumerate}
\end{thm}

\begin{thm}
\cite{p1} The minimal degree sequence for a torus knot of type $(p,\
2p-1),$ $p \ge 2$ denoted by $K_{p,2p-1}$ is given by $(2p-1,\ 2p,\
d)$, where $d$ lies between $2p+1$ and $4p-3$.
\end{thm}

In order to represent a knot-type by a polynomial embedding we
require a suitable {\it knot diagram}. For example for torus knots
of type $(p,q)$ we use its representation as closure of a p-braid
namely $(\sigma_1.\sigma_2\ldots\sigma_{p-1})^q$ and for the
2-bridge knots we use its representation as numerator closure of a
{\it rational tangle}. For a general knot-type there may not be
such systematic nice diagram available.

\begin{defn}
For any two positive integers $p$ and $q$, the $p$-braid $(\sigma_1
\ldots \sigma_{p-1})^q$ is called the {\it toric braid} of type $(p,
q)$.
\end{defn}

Closure of a toric braid gives a torus link of type $(p, q)$. In
particular if $(p, q) = 1$, then we obtain the torus knot of type
$(p,q)$, and it is denoted by $K_{p,q}$.

\begin{defn}
A braid $\beta$ is said to be {\it quasitoric} of type $(p, q)$ if
it can be expressed as $\beta_1 \cdots \beta_q$, where each $\beta_j
= \sigma_1^{e_{j,1}} \cdots \sigma_{p-1}^{e_{j, p-1}}$, with
$e_{j,k}$ is either $1$ or $-1$.
\end{defn}

A quasitoric braid of type $(p, q)$ is a braid obtained from the
standard diagram of the toric $(p, q)$ braid by switching some of
the crossing types.

\begin{thm}(Manturov's Theorem \cite{vom}) Each knot isotopy class can be
obtained as a closure of some quasitoric braid.
\end{thm}

\section{Polynomial Representation of a general knot type}

\begin{thm} \label{t1}
Let $K$ be a knot which is closure of a quasitoric braid obtained
from a toric braid $(\sigma_1 \ \sigma_2 \ \cdots \
\sigma_{p-1})^q$, where $(p, q) =1$, by making $r$ crossing changes.
Then $(2p-1, q + r_0, d)$ is a degree sequence for $K$, where $r_0$
is the least positive integer such that $(2p-1, q + r_0) = 1$ and $d
\le 2q - 1 + 4r$.
\end{thm}

In order to prove this theorem, we first prove the following lemmas.

\begin{lem} \label{l1}
Let $t \mapsto (f(t), g(t))$ represents a regular projection of a
knot $K$. Let $N$ be the number of variations in the nature of the
crossings as we move along the knot. Then there exists a polynomial
$h(t)$ of degree $N$ such that the embedding $t \mapsto (f(t), g(t),
h(t))$ is a representation of $K$.
\end{lem}

\noindent {\bf Proof.} Let $s_1 < s_2 < \ldots < s_N$ be such that
all crossings correspond to parameter values $t \in
\mathbb{R}\setminus \{s_1,s_2,\ldots,s_N\}$ and in any of the open
intervals $(-\infty, s_1)$, $(s_1,s_2)$, $\ldots,\ (s_N,\infty)$ all
the crossings are of the same type (either over or under), also in
successive intervals, the crossings are of opposite type. Now define
\[h(t) = \pm \prod_{i = 1}^{N} (t-s_i).\] It is easy to observe that
$h(t)$ has constant sign on each interval and opposite sign on
consecutive intervals, i.e. it provides an over/under crossing data
for the knot-type. Hence $t \mapsto (f(t), g(t), h(t))$ is a
polynomial representation of $K$.

\begin{lem}\label{l2}
Let $K$ be a knot represented by a polynomial embedding $t \mapsto
(f(t), g(t),h(t))$. Let $K_r$ be a knot obtained from $K$ by making
$r$ crossing changes from over to under or vice versa. Let $N$ be
the degree of $h(t)$ polynomial. Then $K_r$ can be represented by a
polynomial embedding \[t \mapsto (f(t), g(t), h_r(t))\] where
$deg(h_r(t))$ is at most $N + 4r$.
\end{lem}

\noindent {\bf Proof.} It can be easily shown by induction.

{\bf Proof of Theorem \ref{t1}.} Since the regular projection of
$K$, is same as that of $K_{p, q}$, it has $(p-1)q$ real ordinary
double points.

\begin{case} When $(q+1, 2p-1)=1.$ In this case by taking a real
deformation of  the curve $\tilde{C}: (t^{2p-1},t^{q+1})$ with
maximum number of real nodes, which is $(p-1)q$, we obtain a regular
projection of this knot. This deformation exists by Norbart
A'Campo's Theorem \cite{na1}.
\end{case}

\begin{case}
When $(q+1,2p-1)\neq 1.$ Here we choose the least positive integer
$r_0 $ such that $(2p-1, q+ r_0)=1$ and consider the curve
$\tilde{C}:(X(t),Y(t)) = (t^{2p-1},t^{q+r_0})$. The maximum number
of double points in a deformation of this curve is $(p-1)(q+r_0 -1)
= (p-1)q + (r_0 -1)(p-1)$. By a result from {\it real algebraic
geometry }\cite{dp2}, we can choose a real deformation $\tilde{C}:
(X(t), Y(t)) = (f(t), g(t))$ with $deg(f(t)) = 2p-1$, $deg(g(t)) = q
+ r_0$ such that $\tilde{C}$ has $(p-1)q$ real nodes and $(r_0
-1)(p-1)$ imaginary nodes.
\end{case}

Observe that the crossing data for this knot differs from that of
$K_{p, q}$ at r places, by Lemma \ref{l2}, there exists a polynomial
$\tilde{h}(t)$, with degree $d \leq 2q-1 + 4r$, which provides an
over/under crossing data for $K$. Thus we have shown that $(2p-1,
q+r_0, d)$ is a degree sequence for $K$.

\section {The knot $8_{17}$}

Consider the knot $8_{17}$ whose minimal braid representation is
$\sigma_1^2 \sigma_2^{-1} \sigma_1 \sigma_2^{-1} \sigma_1
\sigma_2^{-2}$. Now by using the relation
$\sigma_1\sigma_2\sigma_1=\sigma_2\sigma_1\sigma_1$ we can replace
$\sigma^{-1}\sigma_2$ by
$\sigma_2\sigma_1\sigma_2^{-1}\sigma_1^{-1}$ and obtain an
equivalent braid representation as:

$\sigma_1^{-1}\sigma_2\sigma_1\sigma_2^{-1}\sigma_1^{-1}\sigma_2
\sigma_1\sigma_2^{-1}\sigma_1^{-1}\sigma_2\sigma_1\sigma_2^{-1}
\sigma_1^{-1}\sigma_2.$

This is a quasitoric braid representation for $8_{17}.$ We see that
it is obtained by making crossing changes in the toric braid
$(\sigma_1\sigma_2)^7.$ Thus, a regular projection of $8_{17}$ is
same as a regular projection of a torus knot of type $(3,7).$ To
obtain a regular projection we consider the parametric plane curve:

$(x(t),y(t))= (f(t),g(t))=(t(t^2 - 6.431)(t^2 - 15.91), t(t^2 -
0.18)(t^2 - 2.4899)(t^2 - 17.458)(t^2 - 16.15)(t^2 - 14.8)(t^2 -
11))$.

The projection is shown in figure 4.1 below

\begin{center}
\epsfig{file=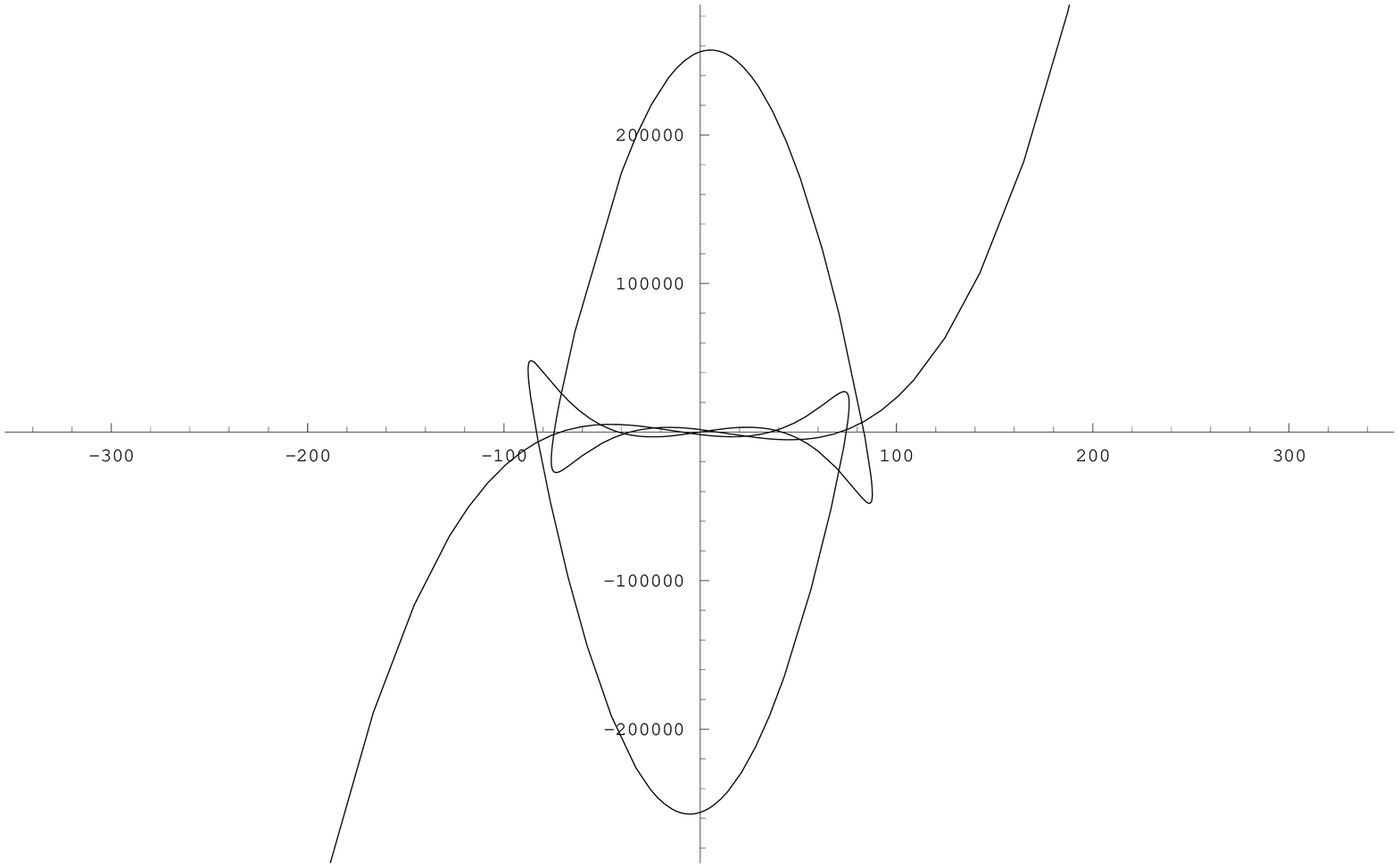,height=12cm,width=14cm}\\
figure 4.1
\end{center}

Observe that $8_{17}$ Knot is obtained from the toric braid
representation of $(3,7)$-torus with 7 crossing changes. From the
quasitoric representation we notice that there are 9 variations in
the signs at the crossing points as we move along the knot. Thus
when we compute the parametric values at the crossing points we
construct the polynomial $h(t)$ using the algorithm in Lemma
\ref{l1}, as $h(t)=(t + 4.138362)(t + 3.86)(t + 2.416735)(t +
1.2)(t)(t - 2.416735)(t - 1.2)(t - 3.86)(t - 4.138362).$ This
polynomial provides us the under/over crossing data for $8_{17}$. A
3 dimensional plot for the curve given by
$(x(t),y(t),z(t))=(f(t),g(t),h(t))$ is shown below in figure 4.2.

\begin{center}
\epsfig{file=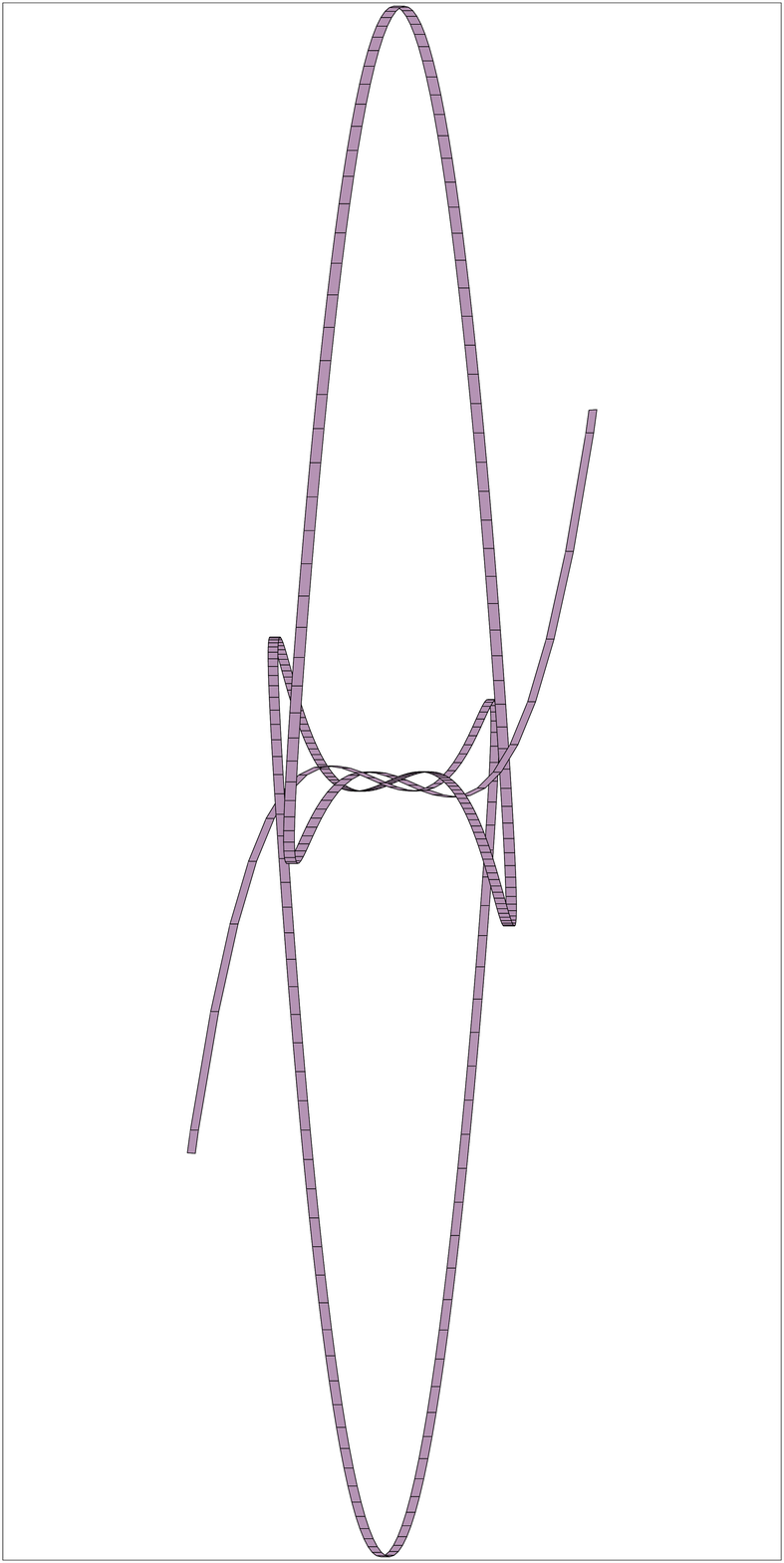,width=7cm}\\
figure 4.2
\end{center}

\section{Polynomial representation of all knots upto 8 crossings}
Here we include  polynomial representations for all knots upto 8
crossings. One can plot their 3 dimensional plots in {\it
Mathematica} or {\it Maple} by giving the ``{\it spacecurve}"
command.
\noindent \textbf{\underline{$3_1$ Knot}:}\

$t \mapsto(t(t - 1) \times (t + 1), t^2(t - 1.15) \times (t +
1.15),(t^2 - 1.056445^2) \times (t^2 - 0.644893^2)t)$

\noindent \textbf{\underline{$4_1$ Knot}:}\

$t \mapsto(t(t - 2) \times (t + 2), (t - 2.1) \times (t +
2.1)t^3,(t^2 - 2.176385^2) \times (t^2 - 1.83588^2) \times (t^2 -
0.8956385^2)t)$

\noindent \textbf{\underline{$5_1$ Knot}:}\

$t \mapsto(t^3 - 4t, (t^2 - 1.2) \times (t^2 - 2.25) \times (t^2 -
3.9) \times (t^2 - 4.85),(t^2 - 2.26311^2) \times (t^2 - 2.116775^2)
\times (t^2 - 1.812575^2) \times (t^2 - 1.2655995^2)t)$

\noindent \textbf{\underline{$5_2$ Knot}:}\

$t \mapsto(t^3 - 17t, t^7 - 0.66t^6 - 29t^5 + 43t^4 + 208t^3 -
680t^2 - 731t,(t + 4.5) \times (t + 4.1) \times (t + 3.2) \times (t
+ 2.3) \times (t + 0.85) \times (t - 0.2) \times (t - 1.75) \times
(t-3.65) \times (t-4.59949))$

\noindent \textbf{\underline{$6_1$ Knot}:}\

$t \mapsto(-16t + t^3, (t + 3.11918) \times (t - 4.519)t(t - 1.95)
\times (t + .9) \times (t + 3.85) \times (t - 1.6125) \times (t +
3.0125) \times (t + 4.35),(t + 4.53430) \times (t + 4.247655) \times
(t + 3.717045) \times (t + 3.23695) \times (t + 2.257325) \times (t
+ 1.22656) \times (t + 0.4656315) \times (t - 0.6227735) \times (t -
1.71873) \times (t - 3.401345) \times (t - 4.51885))$

\noindent \textbf{\underline{$6_2$ Knot}:}\

$t \mapsto((t^2 - 12) \times (t^2 - 11), t(t^2 - 21) \times (t^2 -
7),(t^2 - 4.6573875^2) \times (t^2 - 4.472939^2) \times (t^2 -
3.504525^2) \times (t^2 - 2.318071^2) \times (t^2 - 1.2983325^2)t)$

\noindent \textbf{\underline{$6_3$ Knot}:}\

$t \mapsto(t(t^2 - 16), t^7 - 2.32015t^6 - 37.8493t^5 + 68.303t^4 +
294.038t^3 - 486.111t^2 + 787.942t,(t + 4.385695) \times (t +
4.02568) \times (t + 3.51085) \times (t + 2.11637) \times (t +
1.196606) \times (t + 0.1474538) \times (t - 0.7428872) \times (t -
2.52954) \times (t - 3.810875) \times (t - 4.211535) \times (t -
4.53938))$

\noindent \textbf{\underline{$7_1$ Knot}:}\

$t \mapsto(t^3 - 3t, (t^2 - 1) \times (t^2 - 1.46) \times (t^2 - 2)
\times (t^2 - 3.8) \times (t^2 - 3) \times (t^2 - 4),(t^2 -
1.99354^2) \times (t^2 - 1.962545^2) \times (t^2 - 1.83503^2) \times
(t^2 - 1.56447^2) \times (t^2 - 1.293475^2) \times (t^2 -
1.09503^2)t)$

\noindent \textbf{\underline{$7_2$ Knot}:}\

$t \mapsto(-14t + t^3, -1/ 1000(t + 3.059911918) \times (t - 4.2519)
\times (t + 1.12) \times (t - 2.185) \times (t + .24) \times (t +
4.095) \times (t - 1.19) \times (t - .65125) \times (t + 3.5) \times
(t + 3.935),(t + 4.272945) \times (t + 4.13026) \times (t +
3.841115) \times (t + 3.348575) \times (t + 2.91925) \times (t +
2.126655) \times (t + 1.29503) \times (t + 0.6536465) \times (t -
0.2487545) \times (t - 1.006501) \times (t - 1.72173) \times (t -
3.18092) \times (t - 4.252285))$

\noindent \textbf{\underline{$7_3$ Knot}:}\

$t \mapsto(t^3 - 18\times t, -(t - 4) \times (t + 4) \times (t -
4.472139) \times (t + 4.62139) \times (t^2) \times (t - 3) \times (t
+ 1.86) \times (t^2 - 5.0795) \times (t + 4.975) \times (t + 2.45)
\times (t^2 - 14),(t + 4.854055) \times (t + 4.71158) \times (t +
4.54806) \times (t + 4.23164) \times (t + 2.229356) \times (t +
0.0660095) \times (t - 1.2478465) \times (t - 2.649115) \times (t -
3.458895) \times (t - 3.84593) \times (t - 3.993505) \times (t -
4.25932) \times (t - 4.469185))$

\noindent \textbf{\underline{$7_4$ Knot}:}\

$t \mapsto(t(t^2 - 17), t^2(t^2 - 18) \times (t + 4.7) \times (t^2 -
4.15) \times (t - 4.7),(t^2 - 4.6^2) \times (t^2 - 4.35^2) \times
(t^2 - 4.18^2) \times (t^2 - 9) \times (t^2 - 1.8^2) \times (t^2 -
0.75^2)t)$

\noindent \textbf{\underline{$7_5$ Knot}:}\

$t \mapsto(-22.5t^2 + t^4, -2682.4t + 658t^3 - 46.3t^5 + t^7,(t^2 -
(4.68844)^2) \times (t^2 - (4.42719)^2) \times (t^2 - (4.12207)^2)
\times (t^2 - \ (3.336625)^2) \times (t^2 - (2.579105)^2) \times
(t^2 - (1.3191385)^2)t)$

\noindent \textbf{\underline{$7_6$ Knot}:}\

$t \mapsto((t - 3.25) \times (t + 2.95) \times (t^2 - 18), t(t^2 -
6) \times (t - 3.65) \times (t + 3.45) \times (t^2 - 24),(t +
4.9267) \times (t + 4.83429) \times (t + 4.01544) \times (t +
2.77432) \times (t + 2.11661) \times (t + 1.77356) \times (t +
0.03743) \times (t - 1.894415) \times (t - 2.3901) \times (t -
3.17477) \times (t - 4.3213) \times (t - 4.89905) \times (t -
4.97255))$

\noindent \textbf{\underline{$7_7$ Knot}:}\

$t \mapsto(t(t^2 - 16), t^2(t^2 - 14) \times (t + 4.45) \times (t^2
- 4.85) \times (t - 4.46),(t + 4.46) \times (t + 4.192) \times (t +
3.89) \times (t + 3.355) \times (t + 1.99) \times (t + 0.68) \times
(t - 0.00033) \times (t - 0.663) \times (t - 1.944) \times (t -
3.369) \times (t - 3.888) \times (t - 4.19) \times (t - 4.47))$

\noindent \textbf{\underline{$8_1$ Knot}:}\

$t \mapsto(-13.2t + t^3, -(t + 3.159911918) \times (t - 4.1519)
\times (t + 1.16812) \times (t - 2.1285) \times (t + .24) \times (t
+ 4.095) \times (t - 1.019) \times (t^2 - 2.25) \times (t - .65125)
\times (t + 3.65) \times (t + 3.9035),(t + 4.184905) \times (t +
4.12217) \times (t + 3.989295) \times (t + 3.741405) \times (t +
3.33439) \times (t + 2.84545) \times (t + 2.07536) \times (t +
1.229941) \times (t + 0.510687) \times (t - 0.259909) \times (t -
0.8836655) \times (t - 1.44144) \times (t - 1.885795) \times (t -
3.09575) \times (t - 4.152015))$

\noindent \textbf{\underline{$8_2$ Knot}:}\

$t \mapsto((t^2 - 17.56) \times (t), -(t - 4.09) \times (t - 4.0252)
\times (t + 2.3) \times (t + 4.4954) \times (t + 4.6135499) \times
(t) \times (t + 2.3) \times (t - 3.659028) \times (t - 4.7629)
\times (t + 4.85) \times (t - 1.066352) \times (t + 2.09) \times (t
- 2.19764829) \times (t + 1.1),6.740910^6 + 1.033929 \times 10^7
\times t - 3.1350977\times10^7 \times t^2 + 6.88712\times10^6 \times
t^3 + 1.329509292\times10^7 \times t^4 - 2.647645\times10^6 \times
t^5 - 2.28354726\times10^6 \times t^6 + 350174.87905387 \times t^7 +
200811.8844882988 \times t^8 - 22575.77543179 \times t^9 -
9585.64923687 \times t^{10} + 720.761949 \times t^{11} + 236.8726
\times t^{12} - 9.13684272055 \times t^{13} - 2.378845 \times
t^{14})$

\noindent \textbf{\underline{$8_3$ Knot}:}\

$t \mapsto(t(t^2 - 16), (t^2 - 4.58^2) \times (t^2 - 4.1^2) \times
(t) \times (t^2 - 1.759^2) \times (t^2 - 1.98^2) \times (t^2 - 2^2)
\times (t^2 - 4.13^2),-1.253835\times10^6 \times t +
1.89348\times10^6 \times t^3 - 93478 \times t^5 - 24316.5388 \times
t^7 + 2306.0931 \times t^9 - 54.3551 \times t^{11})$

\noindent \textbf{\underline{$8_4$ Knot}:}\

$t \mapsto(-17.0275 \times t + t^3, -7238.08 \times t +
2156.21\times t^2 + 9252.69 \times t^3 - 2762.63\times t^4 -
2278.64\times t^5 + 686.535 \times t^6 + 154.598 \times t^7 -
47.6818 \times t^8 - 3.15162 \times t^9 + t^{10}),-366161 +
779611.474 \times t + 1.34584\times10^6 \times t^2 -
3.74732\times10^6 \times t^3 + 1.40482\times10^6 \times t^4 +
2.93897\times10^6 \times t^5 - 991572 \times t^6 - 528142 \times t^7
+ 171904.38 \times t^8 + 38975.6 \times t^9 - 12649.82\times t^{10}
- 1290.85\times t^{11} + 425.449\times t^{12} + 15.9088\times t^{13}
- 5.39827\times t^{14})$

\noindent \textbf{\underline{$8_5$ Knot}:}\

$t \mapsto (102.6 \times t - 9.6 \times  t^2 - 22.4125 \times  t^3 +
0.6 \times  t^4 + t^5, 913.915 - 10.1304 \times  t - 1223.62 \times
t^2 + 13.0047 \times  t^3 + 342.138572  \times t^4 - 3.0753 \times
t^5 - 33.43324 \times  t^6 + 0.201  \times t^7 + t^8, -797977.96519
- 1.80074 \times 10^6 \times  t + 5.20615 \times 10^6 \times  t^2 -
3.14955 \times 10^6 \times  t^3 - 1.8510277 \times 10^6 \times  t^4
+ 2.2888687 \times 10^6 \times  t^5 - 2548.54 \times  t^6 -
537208.72308 \times  t^7 + 74673.549 \times  t^8 + 57393.8 \times
t^9 - 11170.2417 \times  t^{10} - 2854.362 \times  t^{11} + 633.4472
 \times t^{12} + 53.41236 \times  t^{13} - 12.658  \times t^{14})$

\noindent \textbf{\underline{$8_6$ Knot}:}\

$t \mapsto(21.56135 - 22.56135 \times t^2 + t^4, 57430.3 \times t^3
- 15828.94409 \times t^5 + 1554.4382 \times t^7 - 65.335787 \times
t^9 + t^{11}, 3.83628\times10^6 \times t - 3.80231\times10^7 \times
t^3 + 2.70379\times10^7 \times t^5 - 7.48485\times10^6 \times t^7 +
1.0666\times10^6 \times t^9 - 85990.4898 \times t^{11} + 3968.12
\times t^{13} - 97.833 \times t^{15} + t^{17})$

\noindent \textbf{\underline{$8_7$ Knot}:}\

$t \mapsto(18.1476 \times t - t^3, 334796 - 355346.8 \times t -
139874 \times t^2 + 174113.2839 \times t^3 + 26617.726 \times t^4 -
28780.456 \times t^5 - 2402.053 \times t^6 + 2153.9814214 \times t^7
+ 99.996212 \times t^8 - 75.2884 \times t^9 - 1.547 \times t^{10} +
t^{11}, 1.135485\times10^7 - 1.06477\times10^7 \times t -
1.61826\times10^7 \times t^2 + 1.597485\times10^7 \times t^3 +
5.23349\times10^6 \times t^4 - 3.97577\times10^6 \times t^5 - 765060
\times t^6 + 438043 \times t^7 + 59464.6 \times t^8 - 25179.9 \times
t^9 - 2548.08 \times t^{10} + 741.202 \times t^{11} + 56.95 \times
t^{12} - 8.8426 \times t^{13} - 0.52 \times t^{14})$

\noindent \textbf{\underline{$8_8$ Knot}:}\

$t \mapsto(6.392595 - 7.04755 \times t - 0.255 \times t^2 + t^3,
261.061 \times t + 63.7029 \times t^2 - 662.543 \times t^3 - 165.758
\times t^4 + 339.24453 \times t^5 + 91.9902 \times t^6 -
60.5558383355 \times t^7 - 17.4661813 \times t^8 + 3.3461 \times t^9
+ t^{10},-11243.145144 + 57682.7791 \times t - 121485 \times t^2 +
129882.5142 \times t^3 - 40537 \times t^4 - 59511.45 \times t^5 +
51887.397 \times t^6 + 5765.9847 \times t^7 - 15323.3 \times t^8 +
706.65613 \times t^9 + 2141.363 \times t^{10} - 149.133 \times
t^{11} - 149.772 \times t^{12} + 6.641753 \times t^{13} + 4.214027
\times t^{14})$

\noindent \textbf{\underline{$8_9$ Knot}:}\

$t \mapsto(t^3-16 \times t, -223891 \times t + 117414 \times t^3 -
21544.8 \times t^5 + 1787.65 \times t^7 - 68.79768 \times t^9 +
t^{11}, -6.2307944\times10^6 \times t + 9.3007326\times10^6 \times
t^3 - 2.61655\times10^6 \times t^5 + 326427 \times t^7 - 21124.3
\times t^9 + 696.354 \times t^{11} - 9.2731 \times t^{13})$

\noindent \textbf{\underline{$8_{10}$ Knot}:}\

$t \mapsto(6.9316 \times t - 2.38335 \times t^2 - 6.01325 \times t^3
+ 0.465 \times t^4 + t^5, -7.404583 + 0.297817 \times t + 13.6921
\times t^2 - 0.307817 \times t^3 - 7.2875 \times t^4 + 0.01 \times
t^5 + t^6, 56.115 - 224.5194 \times t - 466.266 \times t^2 + 609.281
\times t^3 + 1362.664 \times t^4 - 764.397 \times t^5 - 1657.83276
\times t^6 + 480.98677 \times t^7 + 992.291 \times t^8 - 146.86217
\times t^9 - 303.83859 \times t^{10} + 20.8645 \times t^{11} +
45.67642 \times t^{12} - 1.10155 \times t^{13} - 2.676137 \times
t^{14})$

\noindent \textbf{\underline{$8_{11}$ Knot}:}\

$t \mapsto(-t^3+18.1846 \times t, 71592.5 + 53428.5 \times t -
98454.7 \times t^2 - 72895.3 \times t^3 + 25882.5 \times t^4 +
17475.87426 \times t^5 - 2701.076 \times t^6 - 1634.637328 \times
t^7 + 123.0385 \times t^8 + 66.87181 \times t^9 - 2.0261152 \times
t^{10} - t^{11},-2.6439\times10^6 + 9.11332\times10^6 \times t -
886938.328 \times t^2 - 1.29908\times10^7 \times t^3 +
6.06538\times10^6 \times t^4 + 5.47188\times10^6 \times t^5 -
2.228788\times10^6 \times t^6 - 812987.32456 \times t^7 + 305430
\times t^8 + 54586.184 \times t^9 - 19725.3164 \times t^{10} -
1709.51 \times t^{11} + 608.914 \times t^{12} + 20.35213 \times
t^{13} - 7.260873 \times t^{14})$

\noindent \textbf{\underline{$8_{12}$ Knot}:}\

$t \mapsto(t^3-18.1846 \times t, -47870.34 \times t + 65021.98
\times t^3 - 15853.7565 \times t^5 + 1520.73 \times t^7 - 64.260848
\times t^9 + t^{11}, -86291.29115469387 \times t + 542095.213 \times
t^3 + 1.2745\times10^6 \times t^5 - 304536.6 \times t^7 + 26001.8
\times t^9 - 968.456472 \times t^{11} + 13.3503267 \times t^{13})$

\noindent \textbf{\underline{$8_{13}$ Knot}:}\

$t \mapsto(t^3-17.0275, -6878.17 \times t + 1971.48 \times t^2 +
8832.4054 \times t^3 - 2534.74 \times t^4 - 2214.42774 \times t^5 +
638.572 \times t^6 + 157.241 \times t^7 - 45.88093 \times t^8 -
3.38257 \times t^9 + t^{10}, -167328.2 + 548139.2 \times t +
467343.298 \times t^2 - 2.56039\times10^6 \times t^3 +
1.531225328\times10^6 \times t^4 + 2.42578\times10^6 \times t^5 -
902477.197 \times t^6 - 470545 \times t^7 + 155821.49679 \times t^8
+ 37082.845 \times t^9 - 11772.7196 \times t^{10} - 1311.79 \times
t^{11} + 411.0726624818 \times t^{12} + 17.32924 \times t^{13} -
5.444827 \times t^{14})$

\noindent \textbf{\underline{$8_{14}$ Knot}:}\

$t \mapsto(t^3-17.0275 \times t, 6463.5 \times t - 1753.5626153
\times t^2 - 8351.885 \times t^3 + 2267.07033 \times t^4 +
2144.774385 \times t^5 - 583.356 \times t^6 - 160.952 \times t^7 +
43.9894 \times t^8 + 3.63382 \times t^9 - t^{10},15567.9 +
362867.427 \times t - 706636 \times t^2 - 1.1441\times10^6 \times
t^3 + 1.66678 \times10^6 \times t^4 + 1.7645677\times10^6 \times t^5
- 759367.37167 \times t^6 - 394719 \times t^7 + 130227.6 \times t^8
+ 34548.7 \times t^9 - 10301.75 \times t^{10} - 1338.163646615
\times t^{11} + 381.8570624 \times t^{12} + 19.185648 \times t^{13}
- 5.378 \times t^{14})$

\noindent \textbf{\underline{$8_{15}$ Knot}:}\

$t \mapsto(250.125 \times t - 31.75 \times t^3 + t^5, -248.752 +
282.339 \times t^2 - 34.5875 \times t^4 + t^6,-188013.238816 \times
t - 1.43839\times10^7 \times t^3 + 4.80456\times10^6 \times t^5 -
621463.4408 \times t^7 + 39055.2126 \times t^9 - 1194.249 \times
t^{11} + 14.2236338 \times t^{13})$

\noindent \textbf{\underline{$8_{16}$ Knot}:}\

$t \mapsto(t^5-32.5 \times t^3+261 \times t, -293.3125 + 328.875
\times t^2 - 36.5625  \times t^4 + t^6, -39013.8 \times t -
1.66705\times10^7 \times t^3 + 5.35937\times10^6 \times t^5 - 669154
\times t^7 + 40738 \times t^9 - 1212.44 \times t^{11} + 14.1333
\times t^{13})$

\noindent \textbf{\underline{$8_{17}$ Knot}:}\

$t \mapsto(6.5876 \times t - 2.09435 \times t^2 - 5.85825 \times t^3
+ 0.365 \times t^4 + t^5, -7.4045829 + 0.297817 \times t +
13.6920829 \times t^2 - 0.307817  \times t^3 - 7.2875  \times t^4 +
0.01 \times t^5 + t^6, 54.0204 - 151.887 \times t - 457.898 \times
t^2 + 332.938 \times t^3 + 1357.9144 \times t^4 - 358.81513 \times
t^5 - 1684.174 \times t^6 + 204.03671 \times t^7 + 1024.084 \times
t^8 - 53.2616 \times t^9 - 317.8122 \times t^{10} + 5.63785 \times
t^{11} + 48.35914 \times t^{12} - 0.1496988387 \times t^{13} -
2.86492 \times t^{14})$

\noindent \textbf{\underline{$8_{18}$ Knot}:}\

$t \mapsto(t^5-5.5 \times t^3+4.5 \times t, -7.8375 + 14 \times t^2
- 7.35 \times t^4 + t^6, -127.627 \times t + 563.155 \times t^3 -
909.757 \times t^5 + 672.438 \times t^7 - 236.4233 \times t^9 +
38.943 \times t^{11} - 2.4293 \times t^{13})$

\noindent \textbf{\underline{$8_{19}$ Knot}:}\

$t \mapsto(t^5-5.5 \times t^3+4.5 \times t, -7.8375 + 14 \times t^2
- 7.35 \times t^4 + t^6,-10.4337 \times t + 18.5762 \times t^3 -
8.13297 \times t^5 + t^7)$

\noindent \textbf{\underline{$8_{20}$ Knot}:}\

$t \mapsto(-6.5876 \times t + 2.09435 \times t^2 + 5.85825 \times
t^3 - 0.365 \times t^4 - t^5, -7.4045829 + 0.297817 \times t +
13.6920829 \times t^2 - 0.307817  \times t^3 - 7.2875  \times t^4 +
0.01 \times t^5 + t^6, -13.5807 + 53.717 \times t + 94.7779759
\times t^2 - 106.665896 \times t^3 - 102.442 \times t^4 + 76.5135
\times t^5 + 35.0952 \times t^6 - 21.957786 \times t^7 -
3.7440720177 \times t^8 + 2.139 \times t^9)$

\noindent \textbf{\underline{$8_{21}$ Knot}:}\

$t \mapsto(-6.5876 \times t + 2.09435 \times t^2 + 5.85825 \times
t^3 - 0.365 \times t^4 - t^5, -7.4045829 + 0.297817 \times t +
13.6920829 \times t^2 - 0.307817 \times  t^3 - 7.2875 \times  t^4 +
0.01 \times t^5 + t^6, -43.3193 - 39.3746 \times t + 120.193  \times
t^2 + 80.7083 \times t^3 - 122.983 \times t^4 - 54.46748 \times t^5
+ 57.5831 \times t^6 + 14.70042 \times t^7 - 12.4378257 \times t^8 -
1.3658548 \times t^9 + t^{10})$

The 3-d plots of these long knots represented by the above polynomial embeddings, can be downloaded from the following link:

http://www.iitg.ernet.in/prabhakar/myproject.htm

\label{}

\end{document}